\documentclass{amsart}
\usepackage{amsmath, amsfonts, amssymb, amsthm}
\usepackage{enumerate}
\newtheorem{proposition}{Proposition}
\newtheorem{theorem}{Theorem}
\newtheorem{lemma}{Lemma}
\newtheorem{corollary}{Corollary}
\newtheorem*{Conj}{The 3-element Basis Conjecture}
\theoremstyle{definition}
\newtheorem{definition}{Definition}
\newtheorem*{question}{Question}

\usepackage{chngcntr}
\newtheorem{claim}{Claim}
\counterwithin*{claim}{theorem}

\begin{document}

\title
[Subspaces of monotonically normal compacta]
{The basis problem for subspaces of monotonically normal compacta}

\author{Ahmad Farhat}
\address{Mathematical Institute, University of Wroc\l{}aw, pl. Grunwaldzki 2/4, 50-384 Wroc\l{}aw, Poland}
\email{ahmad.farhat@math.uni.wroc.pl}


\begin{abstract}
We prove, assuming Souslin's Hypothesis, that each uncountable subspace of each zero-dimensional monotonically normal compact space contains an uncountable subset of the real line with either the metric, the Sorgenfrey, or the discrete topology.
\end{abstract}

\maketitle

\section*{Introduction}

In \cite{To2}, Todorcevic put forward a deep and outstanding program to inspect basis problems in combinatorial set theory. The blueprint can be described as follows: Given a class of structures $\mathcal{S}$,

\begin{enumerate}
 \item Identify critical members of $\mathcal{S}$, and
 \item Show that a certain list of critical members is complete.
\end{enumerate}

In other words, one wishes to find a \emph{basis} for $\mathcal{S}$. In a specific category, the only part lacking for such a plan to be precise is to have interpretations for the terms ``critical'' and ``complete''. 

In the category of topological spaces, for instance, the morphisms used to interpret completeness may be taken to be embeddings. In this sense, for $\mathcal{S}$ a class in this category, a basis for $\mathcal{S}$ is a subclass $\mathcal{S}_0$ such that each member of $\mathcal{S}$ contains a homeomorphic copy of a member of $\mathcal{S}_0$. Ensuring that the choice of $\mathcal{S}_0$ is minimal (in the sense of cardinality), the members of $\mathcal{S}_0$ are clearly then to be interpreted as the critical objects of $\mathcal{S}$.

Let $\mathcal{S}$ now be the class of uncountable first countable regular spaces. A basis conjecture for $\mathcal{S}$ was first formulated by Gruenhage in \cite{Gr1}. This has come to be called the \emph{3-element basis conjecture for uncountable first countable regular spaces}, as it states that it is consistent that each such space contains a set of reals of cardinality $\aleph_1$ with either the metric, the Sorgenfrey, or the discrete topology. If true, this conjecture would afford a deep understanding of the structure of regular spaces; in fact, the truth of even partial instances of it would have far-reaching consequences and settle long-open problems in set-theoretic topology.

To state the conjecture, let us denote by

\begin{itemize}
 \item $D(\omega_1)$, the discrete space on $\omega_1$,
 \item $B$, a fixed $\aleph_1$-dense subset of the open unit interval, and
 \item $B\times \{0\}$, the corresponding subspace of the split interval (or, double arrow space).
\end{itemize}

\begin{Conj}
 PFA implies that the class of uncountable first countable regular spaces has a three element basis consisting of $D(\omega_1)$, $B$, and $B\times \{0\}$.
\end{Conj}

An excellent discussion of the basis conjecture, specifically for perfectly normal compacta, can be found in \cite{GrM}.
 
T. Eisworth proved in \cite{Ei} that each separable Hausdorff space which is a continuous image of an ordered compactum admits a continuous and at-most 2-to-1 map onto a metric space. Gleaning out ideas from \cite{Fr} and \cite{Gr1}, it follows from this result that PFA implies that the class of uncountable subspaces of monotonically normal compacta has a 3-element basis. In this paper, we prove that actually Souslin's Hypothesis is sufficient to guarantee that each uncountable space having a zero-dimensional monotonically normal compactification contains either an uncountable discrete subspace, an uncountable subspace of the real line, or an uncountable subspace of the Sorgenfrey line. In Section 2 we employ towards this end a direct combinatorial analysis, and in Section 3 we use continuum theoretic techniques to obtain finer structure results for monotonically normal compacta. Finally, we conclude the paper by indicating directions for future research.

This work could be seen as a simple exhibit of the power of M.E. Rudin's latest breakthrough: a proof of Nikiel's Conjecture.

\section{Preliminaries}

The reader is referred to \cite{En} for basic facts stated without proof.
Unless stated otherwise, all sets are assumed to be uncountable, and all spaces Hausdorff. Maps and images refer to continuous functions and their images, respectively. A compactum is a compact space, and a continuum is a connected compactum. A linearly ordered set equipped with the usual open interval topology is called an \emph{ordered space}.

Our departing point for the study of monotonically normal compacta is Nikiel's conjecture \cite{Ni}, proved by M.E. Rudin \cite{Ru}. It states that the class of monotonically normal compacta coincides with the class of continuous images of ordered compacta. We therefore use the terms ``compact monotonically normal" and ``image of a compact ordered space" interchangeably.

Prototypical ordered compacta include the Cantor set in the metric category and the split interval in the separable non-metric category. The split interval is a classical space first described by Alexandroff, and otherwise known by the more common name, the double arrow space. Recall that each compact metric space is an image of the Cantor set. Below we give short proofs of folklore results which establish an analogue of this for the split interval. We note in this context the characterization of separable ordered compacta appearing in \cite{Os1}.

Recall first that a \emph{jump} in a linearly ordered set $K$ is a set $\{a,b\}\in [K]^2$ such that no point of $K$ lies between $a$ and $b$.

\begin{proposition}\label{base_prop}
 Each separable ordered compactum is a continuous image of the split interval.
\end{proposition}

\begin{proof}
 Let $X$ be a separable ordered compactum, and let $\mathcal{S}$ be the set of all maximal convex locally countable open subsets of $X$. Since $X$ is separable, $\mathcal{S}$ is countable, as well as each member of $\mathcal{S}$. Let $X'= X\setminus \bigcup \mathcal{S}$, and notice that $X'$ has no isolated points. Let $G$ be the decomposition of $X'$ into jumps in $X'$ and singletons. The quotient space $X'/G$ is a separable ordered continuum, and hence is homeomorphic to $[0,1]$ (see, e.g., \cite{HY}). It is thus clear that $X'$ is an image of a copy $D$ of the split interval; denote the quotient map by $q$. Notice that for each $S\in\mathcal{S}$ that is neither an upper set nor a lower set of $X$, $\overline{S}=[a_S,b_S]$ for some $\{a_S,b_S\}\in [X]^2$ which is furthermore a jump in $X'$. Moreover, for each $S\in\mathcal{S}$, $\overline{S}$ is a countable compact Hausdorff space, and hence, by the Mazurkiewicz-Sierpi\'{n}ski theorem (\cite{MS}), is homeomorphic to a well-ordered set equipped with the order topology. Each compact countable ordinal is easily seen to be a continuous image of the split interval. To finish the proof, insert a copy of the split interval in each jump in $D$ of the form $\{q^{-1}(a_S), q^{-1}(b_S)\}$, $S\in \mathcal{S}$, append two copies of the split interval to $D$ so that they form an upper set and a lower set of the resulting space, and observe that the resulting space is homeomorphic to the split interval and admits a continuous map onto $X$.
\end{proof}

\begin{proposition}\label{main_prop}
 Each separable monotonically normal compactum is a continuous image of the split interval.
\end{proposition}

\begin{proof}
 Let $X$ be a separable space, $K$ a compact ordered space, and $f:K\to X$ an onto map. Let $D$ be a countable dense subset of $X$, and for each $d\in D$ select a point $d'\in f^{-1}(d)$. $K'=cl(\{d':d\in D\})$ is then a separable ordered compactum which maps onto $X$. By Proposition \ref{base_prop}, $K'$ is an image of the split interval.
\end{proof}

\subsection{The main tool}

A subset $C$ of a linearly ordered set $K$ is said to be \emph{convex} if whenever $a,b,c\in K$ are such that $a,b\in C$ and $a<c<b$, then $c\in C$. An \emph{order component} of a subset $A\subseteq K$ is a subset $C\subseteq K$ which is maximal with respect to the properties ``$C \subseteq A$" and ``$C$ is convex.'' It is easy to see that each order component of each open subset of $K$ is open.

For a set $A$, let $\mathcal{P}(A)$ denote the power set of $A$.

\begin{definition}
 For $f: X\to Y$ a map between two spaces, the function $f^\sharp:\mathcal{P}(X)\to \mathcal{P}(Y)$ is defined by setting $f^\sharp(A)=Y\setminus f(X\setminus A)$ for all $A\subseteq X$.
\end{definition}

\begin{definition}
 Let $K$ be an ordered space, and $X$ and $Y$ be spaces.
\begin{itemize}
 \item A map $f: K\to X$ is said to be \emph{order-light} if for every $x\in X$, each order component of $f^{-1}(x)$ is degenerate.
 \item A surjective map $f: X\to Y$ is said to be \emph{irreducible} if $f(A)\neq Y$ for each proper closed subset $A\subseteq X$. Equivalently, $f$ is irreducible iff for each nonempty open subset $U\subseteq X$, $f^\sharp(U)\neq \emptyset$.
 \item  A map $f: K\to X$ is said to be \emph{reduced} if $f$ is both order-light and irreducible.
\end{itemize}
\end{definition}

\begin{lemma}[\cite{MP}, \cite{Tr}]\label{Tr_Lemma}
 Let $X$ be an image of an ordered compactum. Then there are an ordered compactum $K$ and a reduced map $f:K\to X$.
\end{lemma}

The advantage of having a reduced map $f:K\to X$ is that it is in some sense minimal. First, there are no redundant open sets in the domain of $f$. Moreover, requiring $f$ to be order-light ensures that it does not collapse jumps (as does, for example, the standard map from the split interval onto the unit interval). Thus one would expect a lot of the properties of $X$ to be shared by $K$ if there were a reduced map $f: K\to X$. As for the outline of the more or less straight-forward proof of the lemma: first, Zorn's lemma is used to ensure irreducibility, and second, if $g$ is an irreducible map from an ordered compactum $K'$ onto $X$, then identifying the order components of fibers of $g$ yields an ordered compactum $K$ and an order-light map from $K$ onto $X$.

\section{A combinatorial approach}

A \emph{Souslin line} is an ordered space which is ccc and not separable. By \emph{SH} we denote Souslin's Hypothesis that there are no Souslin lines.

\begin{definition}
 Let $X$ and $Y$ be sets. A function $f:X\to Y$ is said to be 1-to-1 at a point $x\in X$ iff $|f^{-1}f(x)|=1$, or, equivalently, iff $f^\sharp(\{x\})=\{f(x)\}$. If $f$ is 1-to-1 at every point of a subset $A\subseteq X$, $f$ is said to be 1-to-1 on $A$. Obviously, $f$ is 1-to-1 iff $f$ is 1-to-1 on $X$.
\end{definition}

Recall that a \emph{branch} of a tree is a maximal chain in the tree.

\begin{theorem}\label{mainthm}
 Assume SH. Each uncountable monotonically normal compactum contains either an uncountable discrete subspace, an uncountable subspace of real type, or an uncountable subspace of the Sorgenfrey line.
\end{theorem}

\begin{proof}
 Let $X$ be an uncountable monotonically normal compactum, and assume that $X$ contains no uncountable discrete subspace.
 
 Let $X'=X\setminus D$, where $D$ is the set of all locally countable points of $X$. By Lemma 1, there are an ordered compactum $K$ and a reduced map $f:K\to X'$. Since $f$ is irreducible, $K$ contains no locally countable points; in particular, $K$ has no isolated points. We first extend a result of \cite{BT} to the effect that under these conditions, $f$ is 1-to-1 on a dense subset of $K$.
 
 \begin{claim}\label{claim1}
  $f$ is 1-to-1 on an uncountable subset $A\subset K$.
 \end{claim}
 
 \begin{proof}
  Assume $f$ is not 1-to-1 on any uncountable subset of $K$.
 
  If $U$ is a non-empty convex open subset of $K$ with endpoints $a=\inf(U)$ and $b=\sup(U)$, and points $c,d\in U$ are such that 
  $$c< d\text{, }f(c)=f(d)\text{, }f([c,d])\subset f^\sharp(U)\text{, and }f([c,d])\cap \{f(a),f(b)\}=\emptyset\text{,}$$ 
  we say that the interval $(c,d)$ is \emph{internal} in $U$.
  
  Now let $U$ be a convex open subset of $K$ with endpoints $a$ and $b$. To show that there are two points $c,d\in U$ such that $(c,d)$ is internal in $U$ (proved in \cite{BT}), let $U'=(a',b')$ be a convex open subset of $f^{-1}(f^\sharp(U) \setminus \{f(a), f(b)\})$. Since $f$ is reduced, $U'\neq \emptyset$, and since $K$ has no isolated points, $U'$ is uncountable. Similarly, we take $V$ to be a convex open subset of $f^{-1}(f^\sharp(U') \setminus \{f(a'), f(b')\})$, in which case $V\neq \emptyset$ and is uncountable. Since $f$ is not 1-to-1 on any uncountable subset of $U$, there is a point $x\in V$ such that $f$ is not 1-to-1 at $x$. There exist then two distinct points $c$ and $d$ in $f^{-1}(f(x))$ such that $c<d$. Since $f(c)=f(d)=f(x)$ and $x\in V$, $f(c), f(d)\in f^{\sharp}(U')\setminus \{f(a'), f(b')\}$. It follows that $c,d\in U'$. Consequently, since $U'$ was chosen to be convex and $f(c),f(d)\notin \{f(a'),f(b')\}$, $[c,d]\subset U'$. Hence, $f([c,d])\subset f(U')\subset f^\sharp(U) \setminus \{f(a), f(b)\}$. The points $c$ and $d$ are thus as required.
 
  We now construct a Cantor tree $T$ of convex open subsets of $K$ ordered by reverse inclusion. 
  To start, let $L_0 = \{K\}$. Assume now that for some $n\in\omega$, the level $L_n$ has already been constructed, and the elements of $L_n$ form a collection of non-empty pairwise-disjoint convex open subsets of $K$. Let $U$ be an element of $L_n$, and let $a=\inf(U)$ and $b=\sup(U)$. It follows that there exists an interval $U_0=(c,d)$ internal in $U$.
  As $f$ is order-light, $U_0$ is non-empty.
  Since $f([c,d])\cap \{f(a),f(b)\}=\emptyset$, $c\neq a$. Since $a=\inf(U)$, there exists an element $e\in U$ such that $a<e<c$, and so the interval $(a,c)$ is non-empty. Hence, there exists an interval $U_1$ internal in $(a,c)$. 
  It is easy to see that in the construction of $U_1$, we may require that $f(\overline{U_1})\cap \{f(b)\}=\emptyset$, so that $U_1$ is also internal in $U$. 
  Take now $\{U_0, U_1\}$ to be the set of successors of $U$. Define $L_{n+1}$ to be the set of all successors of all elements of $L_n$.
  In this manner, we recursively define all sets $L_n$, $n<\omega$.
  Consider now the tree $T=\bigcup \{L_n: n<\omega\}$, whose levels are the sets $L_n$, $n<\omega$. 
  It is easy to see that $T$ is a Cantor tree.
  
  The tree $T$ obviously has uncountably many branches. 
  For every branch $B$ of $T$, let $U_B = int (\bigcap B)$, the interior of $\bigcap B$. 
  Since $K$ is ccc, for all but at most countably many branches $B$ of $T$, $U_B= \emptyset$. 
  Let $B=\{(a_n,b_n):n < \omega\}$ be a branch of $T$ such that $U_B= \emptyset$. 
  For every $n < \omega$, $f([a_{n+1}, b_{n+1}])\cap \{f(a_n), f(b_n)\} = \emptyset$, which implies that $\bigcap B=\bigcap_{n < \omega}[a_n,b_n]$. As an intersection of a descending sequence of closed intervals, $\bigcap B = [a,b]$ for some $a,b\in K$. 
  Since $int(\bigcap B)=\emptyset$, it follows that $(a,b)=\emptyset$. 
  Moreover, since $f(a_n)=f(b_n)$ for every $n<\omega$, it follows that $f(a)=f(b)$. Since $f$ is order-light, we conclude that $a=b$, and thus that $\bigcap B$ is degenerate. 
  By construction therefore, for any branch $B$ of $T$ that satisfies $U_B= \emptyset$, $f$ is 1-to-1 at $k$, where $\{k\} = \bigcap B$. As $T$ has uncountably many such branches, $f$ is 1-to-1 on an uncountable subset of $K$, and this yields a contradiction.
 \renewcommand{\qedsymbol}{$\diamondsuit$}
 \end{proof}

 \begin{claim}\label{claim}
  $X$ contains an uncountable subspace of the unit interval or of the split interval.
 \end{claim}

 \begin{proof}
  It is readily seen that the irreducibility of $f$ implies that $K$ does not contain an uncountable discrete subspace. $K$ is therefore ccc. Since there are no Souslin lines, $K$ is separable. A separable ordered compactum that has no isolated points is a quotient of the split interval, where equivalence classes are either jumps or one point sets (see the proof of Proposition \ref{base_prop}). It thus follows that the set $A$ from the previous claim contains an uncountable subspace $B$ that can be embedded either in $[0,1]$, or in the split interval. Since $f$ is 1-to-1 on $B$, $B=f^{-1}(f(B))$, so $f\restriction B: B\to f(B)$ is also closed, hence a homeomorphism.
 \renewcommand{\qedsymbol}{$\diamondsuit$}
 \end{proof}
 
 As the top and bottom arrows of the split interval are homeomorphic to the Sorgenfrey line, this concludes the proof.
\end{proof}
  
SH is strictly necessary for Theorem \ref{mainthm} - a left separated subspace of a Souslin line is a counterexample.

\begin{corollary}
 Each uncountable separable monotonically normal compactum contains an uncountable dense subspace which is the disjoint union of a subspace of real type and a subspace of the Sorgenfrey line.
\end{corollary}

\begin{proof}
 Let $X$ be an uncountable separable monotonically normal compactum, and let $X'=X\setminus D$, where $D$ is the set of all locally countable points of $X$. By the proof of Theorem \ref{mainthm}, $X'$ contains an uncountable dense subspace $A$ which can be embedded in a quotient of the split interval where equivalence classes are either jumps or one-point sets. 
 The space $A$ can thus be decomposed into subspaces $B$ and $C$, where $B$ is a copy of a subspace of the real line with the Euclidean topology, and $C$ is a copy of a subspace of the Sorgenfrey line. Clearly, $B\cup D$ contains an uncountable dense subspace $B'$ that can also be embedded into the real line with the Euclidean topology. 
 $B'\cup C$ is thus a dense subspace of $X$ which is as required. We note at the end that either one of the spaces $B'$ and $C$ may be empty.
\end{proof}

We also isolate a simple but useful fact which will be used in the next section, and whose proof can be deduced from that of Theorem \ref{mainthm}.

\begin{proposition}\label{prop}
 Assume SH. Each monotonically normal compactum $X$ satisfies one of the following alternatives:
 \begin{enumerate}[ (a)]
  \item $X$ contains an uncountable discrete subspace.
  \item $X$ is separable.
 \end{enumerate}
\end{proposition}

We note that Proposition \ref{prop} also follows easily from results in \cite{Os}.

\section{More on the structure of separable monotonically normal compacta}

In this section, we refine the study of the structure of monotonically normal compacta, obtaining the canonical basis $\{D(\omega_1), B, B\times \{0\}\}$ for the class of uncountable subspaces of monotonically normal compacta, and specializing this result further.

\begin{definition}
Suppose $X$ is a space, $\langle Y_p \rangle_{p\in X}$ is a sequence of non-empty spaces, and $\langle f_p \rangle_{p\in X}$ is an element of $\prod_{p\in X} C(X\setminus \{p\}, Y_p)$. The \emph{resolution} $R$ of $X$ at each point $p\in X$ into $Y_p$ by $f_p$ is defined by setting $R=\bigcup_{p\in X} \{p\} \times Y_p$ and equipping $R$ with the weakest topology such that each subset of the form 
$$U\otimes V = (\{p\} \times V)\cup \bigcup_{q\in U\cap f_p^{-1}(V)} \{q\}\times Y_q$$
is open in $R$, where $p\in X$, $U$ is an open neighborhood of $p$ in $X$, and $V\subseteq Y_p$ is open.
\end{definition}

Of concern to us here are dendritic resolutions, which we describe after recalling some definitions from continuum theory.

A space $Z$ is called a \emph{dendron} if $Z$ is compact and connected, and for every two distinct points $x$ and $y$ of $Z$, there exists a point $z\in Z$ such that $x$ and $y$ belong to distinct components of $Z\setminus \{z\}$. Any \emph{arc} (i.e., an ordered continuum) is a dendron. Metrizable dendrons are called \emph{dendrites}.

A subset $A$ of a dendron $Z$ is said to be a \emph{strong T-set} in $Z$ if $A$ is non-empty and closed, and each component of $Z\setminus A$ is homeomorphic to the real line.

\begin{definition}
 Let $Z$ be a dendron, $A$ be a non-empty closed subset of $Z$, and $C$ a subset of $A$ such that for each $x\in C$, $Z\setminus \{x\}$ has exactly two components $K_{x,0}$ and $K_{x,1}$. Let $R$ be the resolution of $A$ at each point $p\in A$ into $Y_p$ by $f_p$ where

	$$\left\{\begin{array}{lll}
    p\in A\setminus C & \Rightarrow & Y_p = 1\mbox{, and } f_p \mbox{ is the constant map, and}\\
    p\in C & \Rightarrow & Y_p=2\text{, and } f_p(s) = i \mbox{ iff } s\in K_{p, i}\cap A\mbox{.}
	\end{array}\right.$$

Observe that the exact assignment of $K_{p,0}$ and $K_{p,1}$ for points $p\in C$ does not matter; indeed, if the assignments of $K_{p,0}$ and $K_{p,1}$ for any point $p\in C$ are swapped, the resulting resolutions, as defined above, are homeomorphic. Hence, we may define the \emph{split dendritic resolution} $s(A,C,Z)$ of $A$ by means of $C$ with respect to $Z$ to be the resolution $R$.

Notice that $s(A,C,Z) = (A\setminus C)\cup (C\times \{0,1\})$. The split interval is clearly homeomorphic to $s([0,1],(0,1),[0,1])$. 

 When $Z$ is a copy of $[0,1]$, we denote $s(A,C,Z)$ by $A_C$. The split interval, for instance, is homeomorphic to $I_{(0,1)}$, where $I$ denotes the unit interval $[0,1]$.
\end{definition}

An important structure result in our context is Theorem 3.1 from \cite{NiC}. It states that if $X$ is a zero-dimensional space which is an image of a separable ordered compactum, then there exist a dendrite $Z$ and a strong T-set $A$ in $Z$ such that $X$ is homeomorphic to a split dendritic resolution of $A$ with respect to $Z$. Restating this theorem, each zero-dimensional monotonically normal compactum admits a map $f$ onto a metric space $M$ where each fiber of $f$ consists of at most two points, and where $M$ is a strong T-set in a dendrite. 

\begin{lemma}\label{lem}
 Suppose $A\subset I$ is closed and uncountable, and $C\subset A\cap (0,1)$. $s(A,C,I)$ contains a subspace $L$ homeomorphic to $I_S$ for some subset $S\subset (0,1)$ such that $s(A,C,I) \setminus L$ is countable.
\end{lemma}

\begin{proof}
 Let $\mathcal{U}$ be the collection of all maximal convex open subsets of $I\setminus A$, and let $V=\bigcup\{\overline{U}\setminus U:U\in \mathcal{U}\}$. Notice that by removing a single isolated point from $s(A,C,I)$ for every point of $V\cap (C\cup \{0,1\})$, where $\{0,1\}$ is the set of endpoints of $I$, we obtain a subspace $L$ of $s(A,C,I)$ homeomorphic to $I_S$ for some subset $S\subset (0,1)$. Furthermore, since $\mathcal{U}$ is countable, $V$ is countable, and therefore $s(A,C,I) \setminus L$ is countable.
\end{proof}

Recall that a map $f:X\to Y$ is called \emph{monotone} iff each fiber $f^{-1}(y)$, $y\in Y$, is connected.

\begin{theorem}\label{scdthm}
 Assume SH. Each uncountable subspace of each zero-dimensional monotonically normal compactum contains either an uncountable discrete subspace, an uncountable subspace of real type, or an uncountable subspace of the Sorgenfrey line.
\end{theorem}

\begin{proof}
 Let $X$ be an uncountable subspace of a zero-dimensional monotonically normal compactum $Y$, and assume that $X$ contains no uncountable discrete subspace. It follows that $X$ is ccc, and hence we may assume that $Y$ is also ccc. By Proposition \ref{prop}, $Y$ is separable. 
 By Theorem 3.1 of \cite{NiC}, there exist a dendrite $Z$, a strong T-set $A$ in $Z$, and a subset $C\subset A$ such that $Y$ is homeomorphic to the split dendritic resolution $s(A,C,Z)$. Let $p:s(A,C,Z)\to A$ be the natural projection. 
 Since $p$ is an at most 2-to-1 map, $p(X)$ is uncountable.
 Therefore, by Lemma 2.11 of \cite{NiC}, there exists a copy $I$ of $[0,1]$ in $Z$ such that the subset $p(X)\cap I$ of $Z$ is uncountable. 
 $X$ therefore contains an uncountable subspace $B$ such that $B$ can be embedded in $s(I\cap A, I\cap C, I)$. Since $B$ is uncountable, $I\cap A$ is necessarily uncountable. By Lemma \ref{lem}, we conclude that $X$ contains an uncountable subspace $D$ such that $D$ can be embedded in a space of the form $I_E$, where $E\subset (0,1)$. Clearly then, $D$ must contain an uncountable subspace of real type or an uncountable subspace of the Sorgenfrey line, which finishes the proof.
\end{proof}

Let $\kappa$ be an infinite cardinal. A subset $A\subset \mathbb{R}$ is said to be $\kappa$-\emph{dense} if $|A\cap (a,b)| = \kappa$ for every interval $(a,b)\subset \mathbb{R}$. Theorem \ref{scdthm}, together with Baumgartner's classic result that the Proper Forcing Axiom (PFA) implies that any two $\aleph_1$-dense sets of reals are isomorphic \cite{Ba}, now yield a proof of the 3-element basis conjecture for the class of uncountable subspaces of zero-dimensional monotonically normal compacta. 

We note that SH is not sufficient to prove Theorem \ref{scdthm} in the absence of assumption of zero-dimensionality. 
A subset $H\subset \mathbb{R}^\omega$ is called a \emph{Hurewicz set} if $H$ is uncountable and the intersection of $H$ with any zero-dimensional subset of $\mathbb{R}^\omega$ is countable. As shown by Hurewicz in \cite{Hu}, the existence of a Hurewicz set is equivalent to CH. Assume now CH, and let $H$ be a Hurewicz set and $\overline{H}$ be the closure of $H$ in the Hilbert cube $[0,1]^\omega$. $\overline{H}$ is thus a compact metric space, and hence is monotonically normal.
 Obviously, $H$ contains no uncountable discrete subspace, and no uncountable subspace of the Sorgenfrey line.
 Moreover, by definition, $H$ contains no uncountable subspace of the real line. 
 Thus, recalling that Souslin's Hypothesis is consistent with CH, the uncountable subset $H$ of $\overline{H}$ would serve as a counterexample were Theorem \ref{scdthm} to be stated in the absence of the condition on zero-dimensionality and in the absence of the assumption of $\neg$CH.

In the remainder of this section, we enhance our study of the structure of monotonically normal compacta.

Combining Theorem \ref{scdthm} with Theorem 2.6 of \cite{Gr1}, it follows that under PFA, each non-metrizable separable monotonically normal compactum $X$ contains a subspace of the form $A\times 2$ with the lexicographic order topology, for some uncountable $A\subset [0,1]$. When $X$ is moreover zero-dimensional, we can prove a strengthening of this result without any extra set-theoretic assumptions.

\begin{proposition}
 Each uncountable separable monotonically normal compactum $X$ contains a subspace homeomorphic to $I_S$ for some $S\subset (0,1)$. If $X$ is furthermore non-metrizable and zero-dimensional, then $S$ can be chosen to be uncountable.
\end{proposition}

\begin{proof}
 Let $X$ be an uncountable separable monotonically normal compactum. If $X$ is metrizable, then $X$ contains a copy of the Cantor set $I_{\mathbb{Q}\cap (0,1)}$. If $X$ is non-metrizable and not zero-dimensional, then $X$ contains a non-degenerate component $Y$. By Theorem 1 of \cite{Tr1}, $Y$ is metrizable, and so contains a copy of $I_{\mathbb{Q}\cap (0,1)}$.
 Assume now that $X$ is non-metrizable and zero-dimensional. By Theorem 3.1 of \cite{NiC} there exist a dendrite $Z$, a strong T-set $A$ in $Z$, and $C\subset A$, such that $X$ is homeomorphic to the split dendritic resolution $s(A,C,Z)$. It is easy to see, using Lemma 2.11 of \cite{NiC}, that there is a copy $I$ of $[0,1]$ in $Z$ such that $s(I\cap A, I\cap C, I)$ is not metrizable. $I\cap A$ is obviously uncountable, and so by Lemma \ref{lem}, $s(I\cap A, I\cap C, I)$ contains a subspace homeomorphic to $I_S$ for some subset $S\subset (0,1)$ such that $s(I\cap A, I\cap C, I)\setminus I_S$ is countable. Since $s(I\cap A, I\cap C, I)$ is non-metrizable, $S$ is uncountable, and, obviously, $I_S\subset s(I\cap A, I\cap C, I)\subset X$.
\end{proof} 

Notice that the preceding proposition is also a generalization of the fact that each uncountable compact metrizable space contains a copy of the Cantor set. 

The following is a corollary of Theorem \ref{scdthm}.

\begin{corollary}\label{corcor}
 Assume SH. Each ccc zero-dimensional monotonically normal compactum $X$ satisfies one of the following alternatives:
 \begin{enumerate}[ (a)]
  \item $X$ is metrizable.
  \item $X$ contains an uncountable subspace of the Sorgenfrey line.
 \end{enumerate}
\end{corollary}

\begin{proof}
 Let $X$ be a zero-dimensional ccc monotonically normal compactum, and assume that $X$ is non-metrizable. By Proposition \ref{prop}, $X$ is separable. By Theorem 1 of \cite{Os}, $X$ is perfectly normal. By Lemma 2.1 of \cite{Gr1}, there is an uncountable $Y\subset X$ such that no uncountable subspace of $Y$ is metrizable. Finally, by Theorem \ref{scdthm}, $X$ must contain an uncountable subspace of the Sorgenfrey line.
\end{proof}

The following trichotomy now ensues from Corollary \ref{corcor}. 

\begin{corollary}
 Assume SH. Each zero-dimensional monotonically normal compactum $X$ satisfies one of the following alternatives:
 \begin{enumerate}[ (a)]
  \item $X$ is metrizable.
  \item $X$ contains an uncountable subspace of the Sorgenfrey line.
  \item $X$ contains an uncountable discrete subspace.
 \end{enumerate}
\end{corollary}

The next result provides a condition under which separable monotonically compacta contain a full copy of the split interval. 

\begin{theorem}\label{thmthm}
 Each uncountable separable monotonically normal compactum satisfies one of the following alternatives:
 \begin{enumerate}[ (a)]
  \item $X$ contains an uncountable subspace of real type.
  \item $X$ contains a copy of the split interval.
 \end{enumerate}
\end{theorem}

\begin{proof}
 Let $X$ be an uncountable separable monotonically normal compactum, and assume $X$ contains no uncountable subspace of real type. 
 By Theorem 1 of \cite{Tr1}, each component of $X$ is metrizable. Assume that there exists a non-degenerate component $M$ of $X$. 
 $M$ is then an uncountable compact metric space, and hence contains a copy of the Cantor set. It follows that $X$ contains an uncountable subspace of real type, which is a contradiction. Hence each component of $X$ is degenerate. 
 $X$ is therefore zero-dimensional. 
 By Theorem 3.1 of \cite{NiC} there exist a dendrite $Z$, a strong T-set $A$ in $Z$, and $C\subset A$, such that $X$ is homeomorphic to the split dendritic resolution $s(A,C,Z)$. 
 By Lemma 2.11 of \cite{NiC}, there is a copy $I$ of $[0,1]$ in $Z$ such that $I\cap A$ is uncountable. 
 It is easy to see that the set $(A\setminus C)\cap I$ is at most countable, because otherwise, $X$ would again contain an uncountable subspace of the real line with the Euclidean topology. 
 Therefore, there is a copy of the Cantor set $D$ contained in $C$, so that $X$ contains a copy of $s(D,D,[0,1])$. 
 By removing countably many isolated points from $s(D,D,[0,1])$, we can see that there is a copy of the split interval embedded in $X$.
\end{proof}

\begin{corollary}
 Assume SH. Each uncountable monotonically normal compactum $X$ satisfies one of the following alternatives:
 \begin{enumerate}[ (a)]
  \item $X$ contains an uncountable metrizable subspace.
  \item $X$ contains a copy of the split interval.
 \end{enumerate}
\end{corollary}

\begin{proof}
 Let $X$ be an uncountable monotonically normal compactum, and assume $X$ contains no uncountable metrizable subspace. $X$ therefore contains no uncountable discrete subspace, and so is ccc. By Proposition \ref{prop}, $X$ is separable. By Theorem \ref{thmthm}, $X$ contains a copy of the split interval.
\end{proof}

For any cardinal $\kappa$, $\aleph_0\leq \kappa \leq 2^{\aleph_0}$, there exists a homogeneous separable orderable compactum of weight $\kappa$. This has already been mentioned in the literature; e.g., in \cite{KHa}. Indeed, if $F$ is a subfield of $\mathbb{R}$ of size $\kappa$ and $E=F\cap (0,1)$, then $I_{E}$ is homogeneous and has weight $\kappa$. 

\begin{theorem}\label{Prop2}
 Assume PFA. There is a unique homogeneous separable orderable compactum of weight $\aleph_1$.
\end{theorem}

\begin{proof}
 Let $X$ be a homogeneous separable orderable compactum of weight $\aleph_1$. 
 $X$ contains no isolated points, and so by the proof of Proposition \ref{base_prop}, $X$ is homeomorphic to $A_C$ for some closed subspace $A\subset [0,1]$ and $C\subset A$. Assume without loss of generality that $\inf{A}=0$ and $\sup{A}=1$.
 
 Let $\{U_n:n\in \omega\}$ be the set of all components of $[0,1]\setminus A$. For each $n\in \omega$, $U_n=(a_n,b_n)$, an open interval in $(0,1)$ with endpoints $a_n$ and $b_n$. Let $Y$ be the quotient of $X$ which has as equivalence classes all sets $\{a_n,b_n\}$, $n\in \omega$, and singletons, and denote the quotient map by $g$. $g(X)$ is clearly homeomorphic to some $I_D$, $D\subset [0,1]$. Let $h:I_D\to [0,1]$ be the standard 2-to-1 map.
 
 Since $X$ has no isolated points $C\cap \bigcup\{\{a_n,b_n\}:n\in \omega\}=\emptyset$, and so $g$ is 1-to-1 on $C$. Since $X$ has weight $\aleph_1$, $|C|=\aleph_1$. Since $X$ is homogeneous, $hg(C)$ is $\aleph_1$-dense in $(0,1)$. It is therefore clear that $X$ is homeomorphic to $I_{C\cup E}$, where $E=\{a_n:n\in\omega\}$. Since $C$ is an $\aleph_1$-dense subset of $(0,1)$, so is $C\cup E$.
 
 We have thus proved that any homogeneous separable orderable compactum of weight $\aleph_1$ is homeomorphic to $I_H$ for some $\aleph_1$-dense subset $H\subset (0,1)$. But any two such spaces, $I_J$ and $I_K$ where $J$ and $K$ are $\aleph_1$-dense subsets of $(0,1)$, are homeomorphic. This is because under PFA, any two $\aleph_1$-dense subsets of $(0,1)$ are isomorphic as ordered sets by Baumgartner's result. This finishes the proof. 
\end{proof}

The space from the statement of Theorem \ref{Prop2} is homeomorphic to any space $I_C$ where $C$ is an $\aleph_1$-dense subset of $(0,1)$. We therefore denote it by $I_{\aleph_1}$.

\begin{theorem}
 Assume PFA. Let $X$ be a zero-dimensional ccc monotonically normal compactum of weight at most $\aleph_1$. Then $X$ is the union of a metric space and at most countably many copies of $I_{\aleph_1}$. 
\end{theorem}

\begin{proof}
 Recall that in monotonically normal spaces, cellularity coincides with hereditary cellularity (see, e.g., \cite{Os}). Thus, since $X$ is ccc, $X$ contains no uncountable discrete subspace. Since SH follows from PFA, Proposition \ref{prop} implies that $X$ is separable.

  \begin{claim}
   $X$ admits a decomposition into a countable space $D$ and the elements of an at most countable collection $\mathcal{A}$ of spaces each of which is a split dendritic resolution of $[0,1]$.
  \end{claim}

  \begin{proof}
   By Theorem 3.1 of \cite{NiC}, there exist a dendrite $Z$, a strong T-set $A$ in $Z$, and a subset $C$ of $A$ such that $X$ is homeomorphic to the split dendritic resolution $s(A,C,Z)$. By Lemma 2.11 of \cite{NiC}, if we let $\{x_n:n\in\omega\}$ be a countable dense subset of $A$, then $Z\setminus E_Z \subset \bigcup_{n=1}^\infty [x_0,x_n]$, where $E_Z$ denotes the set of all endpoints of $Z$. Since for any $n\geq 1$, $[x_0,x_n]\cap X$ contains at most countably many isolated points, a consequence of the proof of Theorem 3.1 of \cite{NiC} is that there is a countable family $\mathcal{I}$ of copies of $[0,1]$ in $Z$ such that $X$ is covered by a countable space $D$ and spaces of the form $s(A\cap J, C\cap J, J)$, $J\in \mathcal{I}$.
   \renewcommand{\qedsymbol}{$\diamondsuit$}
  \end{proof}
 
  \begin{claim}
   Each member of $\mathcal{A}$ is either metrizable or a union of a metrizable subspace and a copy of $I_{\aleph_1}$.
  \end{claim}

  \begin{proof}
   Let $L\in \mathcal{A}$, and assume that $L$ is not metrizable. We may assume that $L$ is homeomorphic to $A_C$, where $A\subset [0,1]$ is non-empty, closed, and dense in itself, $C\subset A$ is uncountable, and $L$ contains no isolated points. Let $\mathcal{U}$ be the collection of all maximal convex open subsets $U$ of $[0,1]$ such that $|U\cap C|\leq \aleph_0$. Clearly, $\mathcal{U}$ is at most countable. We may also assume without loss of generality that for every $U\in \mathcal{U}$, $C\cap (cl(U)\setminus U)=\emptyset$.
   
   For each $U\in\mathcal{U}$, $cl(U)$ is homeomorphic to $[0,1]$. As a subspace of the resolution $A_C$, each space $cl(U)_{C\cap U}$, $U\in \mathcal{U}$, is an image of the Cantor set. The space $$E = \bigcup\{cl(U)_{C\cap U}: U\in\mathcal{U}\}$$ is thus metrizable.
   
   Let 
   $$F = A_C\setminus (E \setminus \bigcup\{cl(U)\setminus U: U\in\mathcal{U}\})\text{.}$$
   
   Clearly, $L=E\cup F$. We will show that $F$ is homeomorphic to $I_{\aleph_1}$. Since $L$ is not metrizable, $F$ is not metrizable, and so is non-empty, separable, and has weight $\aleph_1$. Since $E \setminus \bigcup\{cl(U)\setminus U: U\in\mathcal{U}\}$ is open, $F$ is compact, and thus orderable. Since $F$ has no isolated points, $C\setminus \bigcup\mathcal{U}$ is $\aleph_1$-dense in $(0,1)$. This, together with the fact that for all $U\in \mathcal{U}$, $C\cap (cl(U)\setminus U)=\emptyset$, implies that $F$ is homogeneous. By Theorem \ref{Prop2}, $F$ is homeomorphic to $I_{\aleph_1}$.
   \renewcommand{\qedsymbol}{$\diamondsuit$}
  \end{proof}
  
 Combining the two claims concludes the proof.
\end{proof}

\section{Beyond separability and compactness}

The natural point for the continuation of the line of research of this paper is to study the structure of non-separable monotonically normal compacta. In general, the structure of compact spaces which contain uncountable discrete subspaces is little understood. What will make such a study more tractable is the association of monotonically normal compacta with orderable compacta via Nikiel's conjecture. Of specific interest will be the case of first countability. In a subsequent work, we will show how Aronszajn lines, specifically Countryman lines, enter the picture in attempting to find a more nuanced base for the class of uncountable subspaces of first countable monotonically normal compacta. We note that beyond first countability, there seems to be a cut off for this direction of research at countable tightness. The one-point compactification of $D(\omega_1)$, for instance, showcases that the canonical basis from the statement of the 3-element basis conjecture is best possible in this context.

Another line of research we intend to undertake is to find a characterization of those monotonically normal compacta that are premetric compacta of degree at most 2 (that is, compacta which are at most 2-to-1 continuous preimages of metric compacta).

Finally, we end by asking whether it is consistent that there is a 3-element basis for the class of uncountable monotonically normal spaces. As monotone normality is by itself a severe restriction in terms of separation axioms when compared to regularity, this question fits in the natural approach of trying to prove weakenings of the basis conjecture for classes of spaces where conditions stronger than regularity hold (as in \cite{Gr2} and \cite{To3}). Before stating the question, we note that M.E. Rudin, answering a question of S. Purisch, constructed an example of a locally compact monotonically normal space which has no monotonically normal compactification \cite{Ru1}.

\begin{question}
 Is it consistent (with PFA) that there is an uncountable ccc monotonically normal space (not necessarily possessing a monotonically normal compactification) that does not contain a set of reals of cardinality $\aleph_1$ with either the metric or the Sorgenfrey topology?
\end{question}

A positive answer to this question implies that for monotonically normal spaces, we cannot do much without (some sort of) compactness. A negative answer implies of course that the basis conjecture is true for uncountable monotonically normal spaces.

\section*{Acknowledgments}
 The author wishes to thank the referee for corrections to the original draft of the paper and for many helpful remarks.

\bibliographystyle{amsplain}

\end{document}